\documentclass[12pt]{amsart}

\usepackage{amssymb}
\usepackage{amsmath}
\usepackage{amsthm}

\newcommand{\ran}{\mathop{\mathrm{ran}}}

\title{Derivations for a class of matrix function algebras}

\author{Benton L. Duncan}

\address{Department of Mathematics\\
300 Minard Hall\\
North Dakota State University\\
Fargo, ND  58105-5075\\
USA}

\email{benton.duncan@ndsu.edu}



\begin{document}

\theoremstyle{plain}
\newtheorem{thm}{Theorem}
\newtheorem{lem}{Lemma}
\newtheorem{prop}{Proposition}
\newtheorem{cor}{Corollary}

\theoremstyle{definition}
\newtheorem{defn}{Definition}
\newtheorem*{construction}{Construction}
\newtheorem*{example}{Example}

\theoremstyle{remark}
\newtheorem*{question}{Question}
\newtheorem*{acknowledgement}{Acknowledgements}
\newtheorem{rem}{Remark}

\begin{abstract} We study a class of matrix function algebras, here denoted
$\mathcal{T}^{+}(\mathcal{C}_n)$.  We introduce a notion of point
derivations, and classify the point derivations for certain finite
dimensional representations of $\mathcal{T}^{+}(\mathcal{C}_n)$.  We
use point derivations and information about $n \times n$ matrices to
show that every $\mathcal{T}^{+}(\mathcal{C}_n)$-valued derivation
on $\mathcal{T}^{+}(\mathcal{C}_n)$ is inner.
\end{abstract}

\maketitle

Certain matrix function algebras arise in some standard
constructions in the theory of non-selfadjoint operator algebras.
They have been studied as semicrossed product operator algebras, see
in particular \cite{Alaimia:1999} and \cite{Pet-D:1985}.  More
recently this same class of algebras have been realized as directed
graph operator algebras, see \cite[Example 6.5]{Kribs-Power:2003a}.

Directed graph operator algebras have been studied in relation to
the standard commutative example, $A(\mathbb{D})$, the algebra of
holomorphic functions in the unit disk with continuous extensions to
$\mathbb{T}$.  For example the automorphism groups of certain graph
operator algebras have connections to the automorphism group of
$A(\mathbb{D})$, see \cite{Alaimia:1999} when the directed graph is
a cycle and \cite{Dav-Pitts:1998} when the directed graph has a
single vertex.  Another example is the description of the ideal
theory for directed graph operator algebras given in
\cite{Jury-Kribs:2004}.

In this paper we extend descriptions of point derivations of the
disk algebra, see \cite[Section 1.6]{Browder:1969} or
\cite{Browder:1967}, to a notion of point derivations on the
directed graph operator algebras coming from cycles.  We exploit the
structure of these directed graph operator algebras as matrix
function algebras to describe when point derivations occur.

In the first section we define our algebras and establish some
notation.  In the second section we define a notion of point
derivation and prove some preliminary results.  In particular, we
introduce inner point derivations and look for ways to tell whether
a point derivation is inner. In the third section we apply these
results to describe the point derivations of the matrix function
algebras.  We are able to classify when non-inner point derivations
arise depending only on the representation, which is in turn related
to ``points'' in the set of maximal ideals. In the last section we
use the point derivations to show that every derivation is inner.

\section{Background and notation}

We begin by defining the class of matrix function algebras we will
study in this paper.  From directed graphs, they arise as the left
regular representation of the directed graphs with $n$ vertices and
$n$ edges connecting each successive vertex in turn, to form a
single loop, or $n$-cycle. We will use the notation
$\mathcal{T}^{+}(\mathcal{C}_n)$ for these algebras, where $n$ is
the length of the cycle in the algebra.

We can view $\mathcal{T}^{+}(\mathcal{C}_n)$ as a matrix function
algebra of the form \[ \begin{bmatrix} f_{1,1}(z^n) & z f_{1,2}(z^n)
& z^2 f_{1,3}(z^n) & \cdots & z^{n-1}f_{1,n}(z^n) \\ z^{n-1}
f_{2,1}(z^n)
& f_{2,2}(z^n) & z f_{2,3}(z^n) & \cdots & z^{n-2}f_{2,n}(z^n) \\
z^{n-2}f_{3,1}(z^n) & z^{n-1}f_{3,2}(z^n) & f_{3,3}(z^n) & \cdots &
z^{n-3}f_{3,n}(z^n) \\ \vdots & \vdots & \vdots & \ddots & \vdots
\\ zf_{n,1}(z^n) & z^2 f_{n,2}(z^n) & z^{3}f_{n,3}(z^n) &
\cdots & f_{n,n}(z^n) \end{bmatrix} \] where $f_{i,j} \in
A(\mathbb{D})$ for all $ 1\leq i,j \leq n$.

These algebras inherit a matricial norm from the matricial norm on
$A(\mathbb{D})$ as $\mathcal{T}^{+}(\mathcal{C}_n)$ can be viewed as
sitting inside $M_n\otimes A(\mathbb{D})$. We will denote by
$A(z^n)$ the algebra $\{ f(z^n): f \in A(\mathbb{D}) \}.$  Notice
that $A(z^n)$ is a subalgebra of $A(\mathbb{D})$ for all $n$.

In what follows, $A$ will always denote an operator algebra, and by
representation we mean a continuous representation of $A$ as an
algebra of operators acting on a Hilbert space $\mathcal{H}$. We
will denote the elementary matrices with a 1 in the $i$-th diagonal
spot by $e_{ii}$.  The notation $x = [x_{ij}]$ will denote an $n
\times n$ matrix, where the $i$-$j$ entry is $x_{ij}$. We will write
$\ell(i,j)$ for the formula $|i-j|(\mod{n})$.  So that the above
matrix form of $\mathcal{T}^{+}(\mathcal{C}_n)$ can be written as
\[ \{[z^{\ell(i,j)}f_{i,j}(z^n)]: f_{i,j} \in A(\mathbb{D}) \mbox{
for all } 1 \leq i,j \leq n\} . \]

Lastly, for $1 \leq i \leq n-1$ define $Z_i \in
\mathcal{T}^{+}(\mathcal{C}_n)$ as the matrix with $z$ in the
$i$-$(i+1)$ position and zeroes elsewhere. Define $Z_n \in
\mathcal{T}^{+}(\mathcal{C}_n)$ as the matrix with $z$ in the
$n$-$1$ position and zeroes everywhere else.  It is not hard to see
that $\mathcal{T}^{+}(\mathcal{C}_n)$ is generated by the set $\{
e_{ii}, Z_i: 1 \leq i \leq n \}$.  This shorthand will be used later
when dealing with specific matrices.

\section{Noncommutative point derivations}

Some authors take the definition that follows as the definition of
derivation, see Chapter 9 in \cite{Paul:2002} for example.  We use
this notation since we wanted to emphasize the connection between
the derivation and the particular representation. This particular
definition also emphasizes the connections with point derivations
from \cite{Browder:1969} which we exploit in later sections.

\begin{defn} Let $\pi: A \rightarrow B(\mathcal{H})$ be a
representation of $A$.  We say that a continuous linear map $D: A
\rightarrow B(\mathcal{H})$ is a {\em point derivation at $\pi$} if
$D(ab) = D(a) \pi(b) + \pi(a) D(b)$ for all $a,b \in A$.
\end{defn}

Of course the function $D(a) = 0$ is a derivation.  We refer to this
derivation as the trivial, or zero, derivation.  We begin by
identifying a special class of derivations, of which the trivial
derivation is a special case.

\begin{defn} For $\pi:A \rightarrow B(\mathcal{H})$ a
representation of the operator algebra $A$ and for $X \in
B(\mathcal{H})$ we define the function $\delta_X : A \rightarrow
B(\mathcal{H})$ by $\delta_X(a) = \pi(a)X-X\pi(a)$ for all $a \in
A$.\end{defn}

Linearity of $\delta_X$ is obvious.  If we let $\{ a_n \}$ be a
sequence in $A$, then $ \lim (\pi(a)X - X \pi(a)) = \pi(\lim a_n)X -
X \pi(\lim a_n) = \delta_{X}(\lim a_n)$ and hence $\delta_X$ is
continuous. We can also see this by noting that $\| \delta_X \| \leq
2 \|\pi \| \|X \|$.

\begin{prop} If $\pi: A \rightarrow B(\mathcal{H})$ is a
representation and $X \in B(\mathcal{H})$ then the function
$\delta_X$ is a continuous derivation at $\pi$. \end{prop}

\begin{proof}  It remains to show only that $\delta_X$ is a
derivation.  Let $a,b \in A$, then \begin{align*} \delta_X(ab) & =
\pi(ab) X -X \pi(ab) \\ & = \pi(a)\pi(b)X - \pi(a) X \pi(b) + \pi(a)
X \pi(b) - X \pi(a) \pi(b) \\ & =  \pi(a) \delta_X(b) + \delta_X(a)
\pi(b).\end{align*}\end{proof}

\begin{defn} Letting $\pi: A \rightarrow B(\mathcal{H})$ be a
representation of the operator algebra $A$ we say that a derivation
at $\pi$ of the form $\delta_X$ is an {\em inner derivation at
$\pi$}.\end{defn}

If the range of $\pi$ is $\mathbb{C}$ non-trivial inner derivations
do not arise.  However, in a noncommutative setting they often do.

\begin{prop} Assume that $\pi: A \rightarrow B(\mathcal{H})$ is a
representation such that $ \ran \pi$ is not isomorphic to
$\mathbb{C}$.  There exists $X \in B(\mathcal{H})$ such that $
\delta_X \not\equiv 0$. \end{prop}

\begin{proof} It is well known that $B(\mathcal{H})' = \mathbb{C}$.
Since, $\ran \pi$ is not isomorphic to $\mathbb{C}$ there exists $a
\in A$ such that $ \pi(a) \not\in B(\mathcal{H})'$.  In other words,
there is $X \in B(\mathcal{H})$ such that $X\pi(a) \neq \pi(a) X$.
It follows that $ \delta_X$ is nontrivial.\end{proof}

What follows is a theorem that, in certain cases, will allow us to
distinguish the inner derivations from other derivations.

\begin{thm}\label{inner} Let $\pi: A \rightarrow B(\mathbb{C}^{nk})$ be a
representation such that \[ \ran \pi \cong \oplus_{i=1}^n M_k \]
where $k \geq 1$ and $n$ is finite.  A derivation $D$ at $\pi$ is
inner if and only if $D|_{\ker \pi} \equiv 0 $. \end{thm}

\begin{proof}  Assume first, that $D$ is inner at $\pi$.  Then for
$a \in \ker \pi$, $D(a) = \pi(a) X- X \pi(a)$ for some $X \in
B(\mathcal{H})$.  Now $\pi(a) = 0$ and hence $ D(a) = 0$.  As $a$
was an arbitrary element of the kernel the forward direction
follows.

Now suppose that $\ker \pi$ is in the kernel of $D$.  We define a
map $\widehat{D}: \ran \pi \rightarrow B(\mathcal{H})$ by $
\widehat{D} (\pi(a)) = D(a)$. If $\pi(a) = \pi(b)$, then $ a-b \in
\ker \pi$ and hence $D(a-b) = 0$.  It follows, by linearity of $D$,
that $D(a) = D(b)$ and hence the map $\widehat{D}$ is well defined.

Next, \begin{align*} \widehat{D} (\pi(a) \pi(b)) &= \widehat{D}(
\pi(ab)) \\ &= D(ab) \\ &= D(a) \pi(b) + \pi(a) D(b) \\ & =
\widehat{D}(\pi(a)) \pi(b) + \pi(a) \widehat{D}(\pi(b)).
\end{align*} It follows that $\widehat{D}$ defines a derivation on
$\ran \pi$. Further, since $\ran \pi$ is finite dimensional it
follows that $\widehat{D}$ is continuous.

Notice that if $ n=1$ then as $M_k$ is simple every $M_k$ valued
derivation is inner, \cite{Kadison-Ringrose:1997}.  If $n > 1$ we
can use exact sequences of cohomology groups, see
\cite{Johnson:1972} to see that a continuous
$B(\mathbb{C}^{nk})$-valued derivation on \[ \oplus_{i=1}^n M_k \]
is inner.  Hence there is $X \in M_{nk}$ such that
$\widehat{D}(\pi(a)) = \pi(a) X - X \pi(a)$. Since
$\widehat{D}(\pi(a)) = D(a)$ the result now follows.
\end{proof}

The next two propositions give us a short method of checking whether
non-inner derivations can occur at $\pi$.  For an ideal $M$, we
denote by $M^2$ the algebraic ideal generated by elements of the
form $bc$ such that $b, c \in M$.  We will denote the norm closure
of the ideal $M^2$ by $ \overline{M^2}$.

\begin{prop}\label{kernelsquared} If $\ker \pi = \overline{(\ker \pi)^2}$ for a
representation $\pi: A \rightarrow B(\mathcal{H})$, then for a
continuous derivation $D$ at $\pi$, $D|_{\ker \pi} \equiv 0$.
\end{prop}

\begin{proof}  Let $a = bc$ where $b,c \in \ker \pi$.  Then, \begin{align*}
D(a) &= D(bc) \\ &= \pi(b) D(c) + D(b) \pi(c) \\ & = 0 D(c) + D(b) 0
\\ & = 0. \end{align*}  Since $(\ker \pi)^2$ is the ideal generated by
elements of the form $bc$ where $b,c \in \ker \pi$ it follows that
$D|_{(\ker \pi)^2} \equiv 0 $.  Now, continuity of $D$ yields the
result.
\end{proof}

\begin{prop}\label{approximateidentity} If the kernel of the
representation $\pi: A \rightarrow B(\mathcal{H})$ has a bounded
left (right) approximate identity then any continuous derivation $D$
at $\pi$ is identically zero on $\ker \pi$.
\end{prop}

\begin{proof} Let $\{ e_{\lambda} \}$ be a bounded left (right)
approximate identity in $\ker \pi$.  Then for any $f \in \ker \pi$
we know that $ \lim e_{\lambda} f = f$. But notice that $
e_{\lambda} f \in (\ker \pi)^2$ and hence $D(e_{\lambda} f) = 0$ for
all $ \lambda$.  As $D$ is continuous it follows that $D(f)= 0$.  As
$f$ was arbitrary the result follows.\end{proof}

\begin{cor} If $A$ is a $C^*$-algebra and $\pi$ is a
$*$-representation then every derivation at $\pi$ is identically
zero on $\ker \pi$.
\end{cor}

\begin{proof} It is well known that the kernel of a
$*$-representation is a $*$-ideal.  Further every $*$-ideal in a
$C^*$-algebra is a $C^*$-algebra and hence has an approximate
identity.  The result now follows.
\end{proof}

\section{Point derivations on $\mathcal{T}^{+}(\mathcal{C}_n)$}

\begin{defn} For $\lambda \in \overline{\mathbb{D}}$ we define the
representation $\varphi_{\lambda}: \mathcal{T}^{+}(\mathcal{C}_n)
\rightarrow M_n$ by \[ \varphi_{\lambda} (
[z^{\ell(i,j)}f_{i,j}(z^n)]) = [
{\lambda}^{\ell(i,j)}f_{i,j}({\lambda}^n)]. \]
\end{defn}

Notice that for $\lambda \neq 0$ the range of $\varphi_{\lambda}$ is
isomorphic to $M_n$.  It follows that $\ker(\varphi_{\lambda})$ is a
maximal ideal of type $\lambda$, see \cite{Alaimia:1999}.  In
contrast, the range of $\varphi_0$ is the diagonal matrices in
$M_n$. It is not the case that the kernel of $ \varphi_0$ is a
maximal ideal.

The representations of the form $\varphi_{\lambda}$ are enough to
ensure semisimplicity of $\mathcal{T}^{+}(\mathcal{C}_n)$ a well
known result for certain graph operator algebras, see
\cite{Davidson-Katsoulis:2004}, or \cite{Jury-Kribs:2004}, and
semicrossed products, see \cite{Pet:1988}.  We include the proof in
this context for completeness, and since it is not difficult.

\begin{prop} The algebras $\mathcal{T}^{+}(\mathcal{C}_n)$ are semisimple.
\end{prop}

\begin{proof} Let \[ a = [z^{\ell(i,j)}f_{i,j}(z^n)]. \] Assume that $
\varphi_{\lambda} (a) = 0$ for all $0 < |\lambda| < 1$.  Then in
particular, $f_{i,j}(\lambda^n) = 0$ for all $ 0< |\lambda| < 1$.
But since $ f_{i,j}(z^n)$ is analytic in $\mathbb{D}$ and
identically zero on a set containing a limit point in $\mathbb{D}$
then $f_{i,j}(z^n) \equiv 0$ for all $i,j$.  Hence $a = (0)$ and the
result follows. \end{proof}

We now define another important class of representations.

\begin{defn} For $1 \leq i \leq n$ define the representation
$\varphi_{i,0}: \mathcal{T}^{+}(\mathcal{C}_n) \rightarrow
\mathbb{C}$ by \[ \varphi_{i,0} ( [z^{\ell(i,j)}f_{i,j}(z^n)]) =
f_{i,i}(0).
\]\end{defn}

For these representations, the range is $\mathbb{C}$ and hence the
kernels give rise to maximal ideals which, in the notation of
\cite{Alaimia:1999}, are of type $0$.  Notice also, that since the
range is $\mathbb{C}$ there will be no inner derivations at
$\varphi_{0,i}$ for all $1 \leq i \leq n$.  More is actually true.

\begin{prop} For $n \geq 2$, there is no nontrivial point derivation
at $\varphi_{0,i}: \mathcal{T}^{+}(\mathcal{C}_n) \rightarrow M_n$,
where $1 \leq i \leq n$.\end{prop}

\begin{proof} We will prove the result for $i = 1$, the general case
proceeds in a similar fashion. A simple calculation tells us that \[
\ker \varphi_{0,1} = \left\{
\begin{bmatrix} z^nf_{1,1}(z^n) & z
f_{1,2}(z^n) & z^2 f_{1,3}(z^n) & \cdots & z^{n-1}f_{1,n}(z^n) \\
z^{n-1} f_{2,1}(z^n) & f_{2,2}(z^n) & z f_{2,3}(z^n) & \cdots &
z^{n-2}f_{2,n}(z^n) \\ z^{n-2}f_{3,1}(z^n) & z^{n-1}f_{3,2}(z^n) &
f_{3,3}(z^n) & \cdots & z^{n-3}f_{3,n}(z^n) \\ \vdots & \vdots &
\vdots & \ddots & \vdots
\\ zf_{n,1}(z^n) & z^2 f_{n,2}(z^n) & z^{3}f_{n,3}(z^n) &
\cdots & f_{n,n}(z^n) \end{bmatrix} \right\} \] where $ f_{i,j} \in
A(\mathbb{D}) $ for all $i,j$.  Multiplying two general elements of
$\ker \varphi_{0,1}$ together one can verify that $\ker
\varphi_{0,1} = (\ker \varphi_{0,1})^2$.  Using Proposition
\ref{kernelsquared} together with Theorem \ref{inner} we know that
every derivation at $\varphi_{0,1}$ is inner.  But since $\ran
\varphi_{0,i} = \mathbb{C}$ any inner derivation is the zero
derivation.  The result now follows.\end{proof}

Unlike the previous class of representations, the representations
$\varphi_{\lambda}$ give rise to derivations which are inner at $
\varphi_{\lambda}$.  It is the derivations which are not inner at $
\varphi_{\lambda}$ which interest us so we now look at what values
of $ \lambda$ give rise to derivations which are not inner.

\begin{thm} For $|\lambda|<1$ there exist non-inner derivations at
$\varphi_{\lambda}: \mathcal{T}^{+}(\mathcal{C}_n) \rightarrow M_n$.
\end{thm}

\begin{proof} We begin by noticing that the map $F: A(\mathbb{D})
\rightarrow \mathbb{C}$ given by $ F(f) = f'(\lambda)$ is a
continuous linear functional and hence completely continuous
\cite[Corollary 2.2.3]{Eff-Ruan:2000}. In particular we know that
the matricial map $F_{n, \lambda}: M_n \otimes A(\mathbb{D})
\rightarrow M_n$, given by $F_{n, \lambda}([f_{ij}]) =
[f'_{ij}(\lambda)]$ is continuous. Now as
$\mathcal{T}^{+}(\mathcal{C}_n)$ is a subalgebra of $M_n \otimes
A(\mathbb{D})$ we know that  $F_{n, \lambda}$ restricted to
$\mathcal{T}^{+}(\mathcal{C}_n)$ yields a continuous linear map.

We need only show that $F_{n, \lambda}$ is a non-inner derivation at
$\varphi_{\lambda}$.  To see that it is a derivation we will look at
$F_n$ applied to $M_n \otimes A(\mathbb{D})$.  In particular, choose
two elements $f = [f_{ij}],g = [g_{ij}] \in M_n \otimes
A(\mathbb{D})$. Now notice that
\begin{align*} F_{n, \lambda}(fg) & = F_{n, \lambda} \left[
\sum_{j=1}^n f_{ij}g_{jk} \right] \\ & = \sum_{j=1}^n [ f'_{i,j}(
\lambda)g_{jk}(\lambda) +
 f_{ij}(\lambda) g'_{jk}(\lambda)]  \\
& = \left( \sum_{j=1}^n [ f'_{i,j}( \lambda)g_{jk}(\lambda)]\right)
+ \left( \sum_{j=1}^n [f_{ij}(\lambda) g'_{jk}(\lambda)] \right)
\\ & = [f'_{ij}(\lambda)][g_{ij}(\lambda)] +
[f_{ij}(\lambda)][g'_{ij}(\lambda)] \\ & = F_{n,
\lambda}(f)\varphi_{\lambda}(g) + \varphi_{\lambda}(f)F_{n,
\lambda}(g). \end{align*} Restricting to
$\mathcal{T}^{+}(\mathcal{C}_n)$ will not affect the derivation
property and hence $F_{n,\lambda}$ yields a derivation at $
\varphi_{\lambda}$.

Recall the definition of $Z_i$ as the matrix with a z in the
$i$-$(i+1)$ position for $ 1 \leq i \leq n-1$, or the $n$-$1$
position for $i=n$ and zeroes elsewhere.  For $\lambda = 0$ we see
that $F_{n, 0}$ is not inner since $F_{n, 0}(Z_i) \neq 0$ for all
$i$ and yet $\varphi_{0}(Z_i) = 0$, applying Theorem \ref{inner}
verifies the result.

For $0 < |\lambda| <1$ let $f = z- (\lambda)^n$. Notice that $f$ is
an analytic function such that $f'(\lambda^n) \neq 0$ and yet
$f(\lambda^n) = 0$.  Now let $\tilde{f}$ be the element of
$\mathcal{T}^{+}(\mathcal{C}_n)$ given by $[z^{\ell(i,j)}f(z^n)]$.
Notice that $\tilde{f} \in \ker \varphi_{\lambda}$.  However, $F_{n,
\lambda}(\tilde{f}) = [(\lambda)^{\ell(i,j)}n \lambda^{n-1}] \neq
0$. The result now follows as in the case of $\lambda =
0$.\end{proof}

In the special case of point derivations at $\varphi_0$ we are able
to show more.  In analogy with a description of certain homology
groups for the quiver algebras corresponding to a single vertex and
countable edges in \cite{Pop:1998a}, we now show a certain amount of
uniqueness for derivations at $ \varphi_0$.

\begin{prop} Let $D$ be a point derivation at $\varphi_0.$  Then $D$
can be written as $D_0+D_1$ where $D_0$ is inner at $ \varphi_0$,
$D_1$ is a point derivation at $ \varphi_0$, and $D_1(a) \neq 0$
guarantees that $a \in \ker \varphi_0$. Further, $D_1$ is uniquely
determined by the numbers $D_1(Z_i)$ with $ 1 \leq i \leq n$
\end{prop}

\begin{proof} Let $D: \mathcal{T}^{+}(\mathcal{C}_n) \rightarrow M_n$ be a point
derivation at $\varphi_0$.  Notice that $\ran \varphi_0$ is finite
dimensional and hence $\ker \varphi_0$ has a Banach space complement
in $\mathcal{T}^{+}(\mathcal{C}_n)$ which we will denote by $ ( \ker
\varphi_0)^c$. Further, every $ a \in
\mathcal{T}^{+}(\mathcal{C}_n)$ can be written uniquely as $x_a +
y_a$ where $ x_a \in (\ker \varphi_0)^c$ and $ y_a \in \ker
\varphi_0$.  Now there exist $ \lambda _i$ such that $ x_a =
\displaystyle{ \sum_{i=1}^n \lambda_ie_{ii} }$ where $e_{ii}$ is the
elementary matrix with $1$ in the $i$-$i$ position and zero
everywhere else.

We claim that if $ a ,b \in \mathcal{T}^{+}(\mathcal{C}_n)$ then,
with respect to the decomposition above, $ \varphi_0(x_ax_b) = 0 $
if and only if $x_a x_b = 0$. Writing $x_a = \sum_{i=1}^n
\lambda_ie_{ii}$ and $x_b = \sum_{i=1}^n \mu_ie_{ii}$ then,
\[ x_a x_b = \sum_{i=1}^n \lambda_i \mu_i e_{ii} \] and the claim
follows.

Define the map $D_1: \mathcal{T}^{+}(\mathcal{C}_n) \rightarrow M_n$
by letting $D_1(x_a + y_a) = D(y_a)$ with respect to the above
decomposition. We will use the claim in the previous paragraph to
show that $D_1$ is a derivation at $\varphi_{0}$.  Linearity, and
continuity are clear. We need only establish the derivation
property.  Now
\begin{align*} D_1(ab) & = D_1((x_a + y_a)(x_b+y_b)) \\ &= D_1(x_ax_b
+ y_ax_b + x_ay_b + y_ay_b) \\ & = D(y_ax_b) + D(x_ay_b) \\ &=
D(y_a) \varphi_0(x_b) + \varphi_0(y_a)D(x_b) + D(x_a)\varphi_0(y_b)
+ \varphi_0(x_a)D(y_b)
\\ & = D(y_a) \varphi_0(x_b) + \varphi_0(x_a)D(y_b) \\ &= D(y_a) \varphi_0 (x_b +
y_b) + \varphi_0(x_a) D_1(x_b + y_b) \\ & = D_1(x_a + y_a)
\varphi_0(x_b + y_b) + \varphi_0(x_a + y_a) D_1(x_b + y_b) \\ &=
D_1(a) \varphi_0(b) + \varphi_0(a) D_1(b)
\end{align*} and hence $ D_1$ is a derivation at $ \varphi_0$.

Notice that $D_0 = D-D_1$ is an inner derivation since $D-D_1|_{\ker
\varphi_0} = 0$.  It follows that every point derivation at $
\varphi_0$ can be written as an inner derivation and a derivation
which sends $ (\ker \varphi_0)^c$ to zero.

Notice that each derivation of the form $D_1$ is uniquely determined
by the value on $\ker \varphi_0 \setminus \overline{(\ker
\varphi_0)^2}$. A technical calculation shows us that the set $\ker
\varphi_0 \setminus \overline{(\ker \varphi_0)^2}$ is given by \[ \{
\lambda_i Z_i: 1 \leq i \leq n \}. \] The result now follows.
\end{proof}

The previous result relies on a nice decomposition of every element
of $\mathcal{T}^{+}(\mathcal{C}_n)$ which is invariant under
derivations. Although we expect a similar result for the point
derivations at $ \varphi_{\lambda}$ for all $0 < |\lambda|< 1$ we
have not been able to prove such a result.

\begin{thm}\label{inneratt} For $\lambda \in \mathbb{T}$ every
derivation at $\varphi_{\lambda}: A \rightarrow M_n$ is
inner.\end{thm}

\begin{proof}  We will show that $\ker \varphi_{\lambda}$ has a
bounded approximate identity and then apply Proposition
\ref{approximateidentity} and Theorem \ref{inner}.

We let $ \pi_n: A(\mathbb{D}) \rightarrow A(\mathbb{D})$ be the
contractive homomorphism induced by sending $ z \mapsto z^n$. Denote
the range of this map by $A(z^n)$ which matches our previous
definition of $A(z^n)$.  Further $\pi_n$ is a contractive
isomorphism onto $A(z^n)$.  (We are not making any claims about
contractivity of the reverse map). Notice that \[ \pi_n (\{ f \in
A(\mathbb{D}): f(\lambda) = 0 \}) \subseteq \{ g \in A(z^n):
g(\lambda^{\frac{1}{n}}) = 0 \}.\] Further, since $ | \lambda | = 1$
we know that there is a uniformly bounded net, see \cite[Section
1.6]{Browder:1969}, \[ \{ f_{\iota} \} \subseteq \{ f \in
A(\mathbb{D}): f(\lambda) = 0 \}\] such that $ f_{\iota} g
\rightarrow g$ for all
\[ g \in \{ f \in A(\mathbb{D}): f(\lambda) = 0 \}.\]  Notice that $
\| \pi_n(f_{\iota}) \| \leq \| f_{\iota} \|$ and hence $ \{
\pi_n(f_{\iota}) \} $ is a bounded net in $\{ g \in A(z^n):
g(\lambda^{\frac{1}{n}}) = 0 \}$.  Now if $g(\lambda^{\frac{1}{n}})
= 0$ and $g \in A(z^n)$ then \[ h = \pi_n^{-1}(g) \in \{ f \in
A(\mathbb{D}): f(\lambda) = 0 \}.\] It follows that $ f_{\iota} h
\rightarrow h$.  Now $ \pi_n(f_{\iota} h) \rightarrow g$ and hence
the ideal \[ \{ g \in A(z^n): g(\lambda^{\frac{1}{n}}) = 0 \}\] has
a bounded approximate identity.

We define the net $\{ F_{\iota} \}$ to be the diagonal matrices with
$\pi_n(f_{\iota})$ along the diagonals.  Now $ \{ F_{\iota} \}$ is a
bounded net as $\{ f_{\iota} \}$ is.  Further, $F_{\iota} \in \ker
\varphi_{\lambda}$ for all $\iota$.  It is not hard to see that $\{
F_{\iota} \}$ is an approximate identity in $\ker
\varphi_{\lambda}$.\end{proof}

We will use this theorem to show the main result in this paper, that
every $\mathcal{T}^{+}(\mathcal{C}_n)$-valued derivation is inner.

\section{Derivations on $\mathcal{T}^{+}(\mathcal{C}_n)$}

We begin with an elementary lemma relating derivations and point
derivations.

\begin{lem} Let $D: A \rightarrow A$ be a continuous derivation on
the operator algebra $A$.  For a representation $\pi: A \rightarrow
B(\mathcal{H})$, the map $\pi \circ D : A \rightarrow
B(\mathcal{H})$ is a continuous derivation at $\pi$.\end{lem}

\begin{proof} Since $ \pi \circ D$ is a composition of continuous
linear maps, it follows that $ \pi \circ D$ is a continuous linear
map.  Now let $a, b \in A$.  Then \begin{align*} \pi \circ D (ab) &
= \pi ( D(a)b + aD(b)) \\ & = \pi(D(a))\pi(b) + \pi(a) \pi(D(b)) \\
&= \pi \circ D (a) \pi(b) + \pi(a) \pi \circ D (b). \end{align*} It
follows that $\pi \circ D$ is a derivation at $\pi$.\end{proof}

\begin{defn}  Let $\pi: A \rightarrow B(\mathcal{H})$ be a
representation and $D: A \rightarrow A$ be a continuous derivation.
We say that $D$ is {\em locally inner at $\pi$} if $\pi \circ D$ is
inner at $\pi$.\end{defn}

We are now in a position to tackle the main theorem of this paper.
Showing that every $\mathcal{T}^{+}(\mathcal{C}_n)$-valued
derivation on $\mathcal{T}^{+}(\mathcal{C}_n)$ is inner will be a
simple corollary.

\begin{thm}\label{locallyinner} Let $D: \mathcal{T}^{+}(\mathcal{C}_n)
\rightarrow \mathcal{T}^{+}(\mathcal{C}_n)$ be a continuous
derivation which is locally inner at $\varphi_{\lambda}$ for all
$\lambda \in \mathbb{T}$, then $D$ is inner.\end{thm}

\begin{proof} For $\lambda \in \mathbb{T}$ let
$D_{\lambda}:= \varphi_{\lambda} \circ D$.  Then, by hypothesis,
$D_{\lambda}(a)$ can be written as $ X_{\lambda}
\varphi_{\lambda}(a) - \varphi_{\lambda}(a) X_{\lambda}$ for all $a
\in \mathcal{T}^{+}(\mathcal{C}_n)$.

Notice that as $e_{ii} \in \mathcal{T}^{+}(\mathcal{C}_n)$ it
follows that $D(e_{ii}) \in \mathcal{T}^{+}(\mathcal{C}_n)$.  By
hypothesis $D_{\lambda}(e_{ii}) = X_{\lambda}e_{ii} -
e_{ii}X_{\lambda}$.  Now $X_{\lambda}e_{ii} - e_{ii}X_{\lambda}$ is
the matrix with $(-X_{\lambda})_{ij} $ in the $i$-$j$ position, for
$i\neq 0$, $(X_{\lambda})_{ji}$ in the $j$-$i$ position for $i \neq
j$ and 0 elsewhere.  In particular, the matrix $Y$ with $0$ on the
diagonal such that $\varphi_{\lambda}(Y_{ij}) = (X_{\lambda})_{ij}$
off the diagonal for all $\lambda \in \mathbb{T}$ is an element of $
\mathcal{T}^{+}(\mathcal{C}_n)$.

Recall the definition of the matrices $Z_i$. Now $D_{\lambda}(Z_i) =
\lambda(X_{\lambda} e_{i,i+1} - e_{i,i+1}X_{\lambda})$ for all $
\lambda \in \mathbb{T}$ where we define $e_{n, n+1}$ to mean
$e_{n,1}$. Now the $i$-$(i+1)$ entry of $D_{\lambda}(Z_i)$ is $
\lambda (-X_{i+1,i+1}(\lambda) + X_{i,i}(\lambda))$, and hence for
all $i$, $X_{i+1,i+1}(\lambda) - X_{i,i}(\lambda)$ defines an
element of $\mathcal{T}^{+}(\mathcal{C}_n)$, call it $X_i$. Now
define a diagonal matrix $X'$ by letting the $i$-$i$ entry be $X_i -
X_1$.

Define $X = Y + X'$ which is in $\mathcal{T}^{+}(\mathcal{C}_n)$.
Further, $D_{\lambda}(e_{ii}) = \varphi_{\lambda}(e_iiX-Xe_ii)$ and
$D_{\lambda}(Z_i) = \varphi_{\lambda}(Z_iX-XZ_i)$.  It follows that
for any $ a \in \mathcal{T}^{+}(\mathcal{C}_n)$, $ D_{\lambda}(a) =
\varphi_{\lambda}(aX-Xa)$ for all $ \lambda \in \mathbb{T}$. Now,
every element of $a \in \mathcal{T}^{+}(\mathcal{C}_n)$ is uniquely
determined by the values of $ \varphi_{\lambda}(a) $ for $ \lambda
\in \mathbb{T}$. It follows that $D(a) = aX-Xa$ and the result is
established.
\end{proof}

\begin{cor}\label{graphinner} Every derivation $D: \mathcal{T}^{+}(\mathcal{C}_n)
\rightarrow \mathcal{T}^{+}(\mathcal{C}_n)$ is inner. \end{cor}

\begin{proof}  First, since $\mathcal{T}^{+}(\mathcal{C}_n)$ is semisimple we
know, \cite{Johnson-Sinclair:1968} that every derivation is
automatically continuous.  We know from Theorem \ref{inneratt} that
every continuous derivation at $ \varphi_{\lambda}$ is inner. In
particular, $D_{\lambda}$ is continuous and locally inner on
$\mathbb{T}$.
\end{proof}

We remark that, as $A(\mathbb{D})$ is a special case of
$\mathcal{T}^{+}(\mathcal{C}_n)$, the above proof is an alternate
approach to the fact there are no nontrivial derivations on
$A(\mathbb{D})$.  It would be interesting to know if the above
result can be extended to $\mathcal{L}_{\mathcal{C}_n}$ which is the
matrix function algebra as $\mathcal{T}^{+}(\mathcal{C}_n)$ with
$A(z^n)$ replaced by $H^{\infty}(z^n)$.


\begin{thebibliography}{00}

\bibitem{Alaimia:1999} M.\ Alaimia, Automorphisms of some Banach
algebras of analytic functions, {\em Linear Algebra Appl.\ } {\bf
298} (1999) 87-97.

\bibitem{Browder:1967} A.\ Browder, Point derivations on function
algebras, {\em J.\ Funct.\ Anal.\ } {\bf 1} (1967) 22-27.

\bibitem{Browder:1969} A.\ Browder, ``Introduction to Function
Algebras,'' W.\ A.\ Benjamin, New York, 1969

\bibitem{Davidson-Katsoulis:2004} K.\ Davidson and E.\ Katsoulis,
Nest representations of directed graph algebras, preprint.

\bibitem{Dav-Pitts:1998} K.\ Davidson and D.\ Pitts,
The algebraic structure of non-commutative analytic {T}oeplitz
algebras, {\em Math.\ Ann.\ } {\bf 311 (2)} (1998), 275-303.

\bibitem{Pet-D:1985} L.\ DeAlba and J.\ Peters, Classification of
semicrossed products of finite-dimensional {$C\sp \ast$}-algebras,
{\em Proc.\ Amer.\ Math.\ Soc.\ } {\bf 95(4) } (1985), 557-564.

\bibitem{Eff-Ruan:2000} E.\ Effros and Z.\-J.\ Ruan ``Operator
Spaces'' Clarendon Press, Oxford, 2000.

\bibitem{Johnson:1972} B.\ Johnson ``Cohomology in Banach
Algebras'' American Mathematical Society, Providence, 1972.

\bibitem{Johnson-Sinclair:1968} B.\ Johnson and A.\ Sinclair,
Continuity of derivations and a problem of {K}aplansky, {\em Amer.\
J.\ Math.\ }{\bf 90} (1968), 1067-1073.

\bibitem{Jury-Kribs:2004} M.\ Jury
and D.\ Kribs, Ideal structure in free semigroupoid algebras from
directed graphs, {\em J.\ Operator Theory} to appear.

\bibitem{Kadison-Ringrose:1997} R.\ Kadison J.\ Ringrose,
``Fundamentals of the Theory of Operator Algebras'', American
Mathematical Society, Providence, 1997.

\bibitem{Kribs-Power:2003a} D.\ Kribs and S.\ Power, Free
semigroupoid algebras, {\em J.\ Ramanujan Math.\ Soc.\ } {\bf 19}
(2004), 75-114.

\bibitem{Paul:2002} V.\ Paulsen, ``Completely Bounded Maps and
Operator Algebras'' Cambridge University Press, Cambridge, 2002.

\bibitem{Pet:1988} J.\ Peters, The ideal structure of certain
nonselfadjoint operator algebras, {\em Trans.\ Amer.\ Math.\ Soc.\ }
{\bf 305(1)} (1988), 333-352.

\bibitem{Pop:1998a} G.\ Popescu, Noncommutative joint dilations
and free product operator algebras, {\em Pacific J.\ Math.\ } {\bf
186(1)} (1998), 111-140.





\end{thebibliography}
\end{document}